\documentclass[a4paper,10pt]{amsart}
\usepackage[utf8]{inputenc}
\usepackage[foot]{amsaddr}

\usepackage[english]{babel}
\usepackage{amsmath,amssymb,verbatim,mathrsfs,latexsym,paralist}
\usepackage{graphicx}
\usepackage{epstopdf}
\usepackage{indentfirst}
\usepackage{amsthm}
\usepackage{hyperref}
\allowdisplaybreaks[4]

\newcommand{\Z}{\mathbb{Z}}
\newcommand{\C}{\mathbb{C}}

\newcommand{\mf}{\mathfrak}
\newcommand{\g}{\mf{g}}
\newcommand{\h}{\mf{h}}

\newcommand{\rk}{\mathrm{rk}}

\numberwithin{equation}{section}
\newtheorem{theorem}{Theorem}[section]
\newtheorem{proposition}[theorem]{Proposition}
\newtheorem{conj}[theorem]{Conjecture}
\newtheorem{lemma}[theorem]{Lemma}

\theoremstyle{remark}
 
\theoremstyle{remark}
\newtheorem{rmk}[theorem]{Remark}
   
\title{ Reeder's Conjecture for Even Orthogonal Lie algebras.}
\author{Sabino Di Trani}
\address{ Dipartimento di Matematica ``Guido Castelnuovo'', Sapienza - Universit\`a di Roma.
}
\email{sabino.ditrani@uniroma1.it }
\address{\emph{The author has been partially supported by GNSAGA - INDAM group.}}
\address{ORCID id: https://orcid.org/0000-0002-6651-558X}

\begin{document}

\maketitle

\textbf{Abstract:}
In the paper we complete a case by case proof of Reeder's Conjecture started in our previous work, proving the conjecture for simple Lie algebras of type $D$ and for the exceptional cases.\\

\textbf{Mathematics Subject Classification (2010):} 17B10, 17B20, 17B22.\\

 \textbf{Acknowledgments:} I am grateful to Professor Paolo Papi thank to whom I focused my attention on Reeder's Conjecture during my doctoral studies and who helped me with the computations for the exceptional cases.
 Moreover I would like to thank Professor Kyo Nishiyama for sharing the  computations of \cite{GnS} for types $E$ and $F_4$.
Finally I would like to thank Professor Martina Lanini, Professor Andrea Maffei, Francesco Ferrante and Alessandro Iraci for many useful discussions about computational issues and optimization strategies that we used in the $E_8$ and $F_4$ cases. \\

\textbf{Data Availability:} Data sharing not applicable to this article.

\section{Introduction}
The structure of $\g$- representation of the exterior algebra $\Lambda \g$, $\g$ a simple Lie algebra over $\C$, has been extensively studied in the last fifty years. Despite its finite dimensionality, an uniform description of irreducible components of $\Lambda \g$ is known only in a conjectural form. 

In '97 Reeder proved that for the irreducible representations $V_\lambda$ indexed by certain dominant weights $\lambda$ called \emph{small} (i.e., such that $\lambda$ is in the root lattice and $2 \alpha$  is not smaller than $\lambda$ in the dominant order for all positive roots $\alpha$), the multiplicity of $V_\lambda$ in  $\Lambda \g$ equals 
$2^{\rk \g } m_\lambda^0$,  where $m_\lambda^0$ is the dimension of the zero weight space in $V_\lambda$. 
Moreover he conjectured that the problem of determining the graded multiplicities of small representations in the exterior algebra can be reduced  to a problem involving Weyl group representations on the zero weight space $V_{\lambda}^0$.   

The setup for the conjecture is the following: 
let $\g$ be a simple Lie algebra over $\C$, fix a Cartan subalgebra $\h$ and let $\Phi$ be the associated root system, with Weyl group $W$. 
We denote by $( \_\, ,\,\_ )$ the $W$-invariant positive-definite inner product on $\h^*$ induced by the Killing form and with $\alpha^\vee$ the coroot associated to $\alpha$. 
We choose a  set of positive roots $\Phi^+$ associated to a simple system $\Delta$. Let $\rho$ be the corresponding Weyl vector and $\theta$ (resp. $\theta_s$) the highest root (resp. the highest short root). We will denote by $\Pi^+$ the set of dominant weights, $\omega_i$ will be the $i$-th fundamental weight.
Let $H$ (resp $H^h$) be the space of $W$-harmonic polynomials (resp. of degree $h$) on $\h$, i.e. polynomials annihilated by constant coefficients $W$-invariant differential operators with positive degree. 
\begin{conj}[Reeder]\label{RedC}
Consider the two polynomials 
 \begin{equation*}
 P(V_\lambda,\bigwedge \g, u)=\sum_{n\geq 0}\dim Hom_\g(V_\lambda,\bigwedge^n \g)u^n,
\end{equation*}
\begin{equation*}
 P_W(V^0_\lambda,q,y)=\sum_{n\geq 0}\dim Hom_W(V^0_\lambda,\bigwedge^k\h\otimes H^h) q^hy^k.
\end{equation*}
If $V_\lambda$ is a small representation, then the following equality holds:
\begin{equation*}
 P(V_\lambda,\bigwedge \g, q)= P_W(V^0_\lambda,q^2,q)
\end{equation*}
\end{conj}
This conjecture was implicitly proved in type $A_n$ comparing the ``First Layer Formulae'' of Stembridge \cite{St1} with the results due to Kirillov, Pak \cite{KP} and Molchanov \cite{M}. Moreover some general formulae for graded multiplicities are proved when $V_\lambda$ is the adjoint or the little adjoint module in \cite{Baz}, \cite{DCPP},  \cite{DCMPP} and \cite{Stembridge}.
In these cases the description of the zero weight space is quite simple and the proof of the Reeder's Conjecture is completely straightforward using the formulae proved in \cite{GnS}.
In our previous paper \cite{SDT} we proved the conjecture for Lie algebras of type $B$ and $C$ using the tools introduced by Stembridge in \cite{Stembridge}.
The aim of this note is to complete a case by case proof of the Reeder's Conjecture for the classical Lie algebras showing that the conjecture holds for even orthogonal Lie algebras. 
Furthermore, we checked using SageMath that the conjecture is true in the exceptional cases. 

The first sections of the paper are dedicated to explain our tools and to make more explicit the ``Weyl group part'' of the conjecture.
Section 4 is devoted to prove Conjecture \ref{RedC} for even orthogonal Lie algebras. We use a mixed strategy with respect to what we have done in our previous paper: we use the combinatorics of weights and the action of the Weyl group to find nice closed expressions for the coefficients of Stembridge's minuscule recurrences and again we conclude using an inductive reasoning. 
Section 5 contains some technical results about a suitable reduction of the Stembridge's recurrences. 
Finally, in Section 6 we give an overview about small representations in the exceptional cases and describe the computational approach that we followed to check the conjecture. 
%
\section{Stembridge's Recurrences}
Our work make an extensive use of results exposed in \cite{Stembridge} about the the coefficients $C_\mu (q,t)$ in the character expansion of Macdonald kernels. The reason for our interest in these functions is that the evaluation $C_\mu (-q,q^2)$ gives the polynomial of graded multiplicities of $V_\mu$ in the exterior algebra.

Let  $\Delta(q,t)$ denote the Macdonald kernel and define $C_\mu (q,t)\in \C[q^{\pm 1},t^{\pm 1}], \, \mu \in \Pi^+$ by the relation $\Delta(q,t)=\sum_{\mu\in  \Pi^+}  C_\mu (q,t) \chi(\mu)$. 
It is possible to extend the definition of $C_\mu (q,t)$ to any weight $\mu$ setting 
 \begin{equation}\label{RedClambda}
C_\mu(q,t) =
\left\{
	\begin{array}{lll}
	0 & \mbox{if } \mu + \rho \mbox{ is not regular}, \\
    (-1)^{l(\sigma)}C_\lambda(q,t)  & \mbox{if } \sigma(\mu + \rho) = \lambda+\rho \; , \lambda \in \Pi^+, \sigma \in W .\\
	\end{array}
\right.
\end{equation}
We will say that, if there exists $\sigma$ such that $\sigma(\mu + \rho)=\lambda + \rho$, the weight $\mu$ is \emph{conjugated} to $\lambda$ and we will write $\mu + \rho \sim \lambda + \rho$.
This rational functions $C_\mu(q,t)$ satisfy some recurrences. The problem of their explicit computation reduces to solving a linear system of equations with coefficients in $ \C[q^{\pm 1},t^{\pm 1}]$.

We recall that a weight (resp. coweight) $\omega$ is said to be \emph{minuscule} if $(\omega, \alpha^\vee ) \in \{0, \pm1\} $ (resp. $(\omega, \alpha)$) for all positive roots.
Fix a dominant weight $\lambda$. If $\omega $ is a minuscule coweight, then the following relation holds (see \cite{Stembridge}, formula (5.14)):
 \begin{equation}\label{ricorsionemin}
  \sum_{i=1}^k C_{w_i \lambda}(q,t) \left(\sum_{\psi \in O_\omega} \left(t^{-(\rho, w_i \psi)}-q^{(\lambda, \omega)}t^{(\rho, w_i\psi)}\right)\right)=0.
 \end{equation} 
 Here $w_1, \dots , w_k$ are minimal coset representatives of $W/W_\lambda$, where $W_\lambda$ is the stabilizer of $\lambda$, and $O_\omega$ is the orbit $W_\lambda \cdot \omega$. We are going to call this recursive relation the \emph{minuscule recurrence}.
The $C_\mu$ appearing in  (\ref{ricorsionemin}) are not necessary in their reduced form (i.e. the weight $\mu$ is not necessary dominant), but the reduced form can be always achieved according to the Definition \ref{RedClambda}.
Considering only the reduced forms, Stembridge proves that the $C_\mu(q,t)$ appearing in (\ref{ricorsionemin}) are indexed only by weights $\mu$ smaller or equal to $\lambda$ in the dominant order. We recall that if $\lambda$ is small, a weight $\mu \leq \lambda$ is again small.
Consequently our strategy becomes very clear:
we determine closed formulae for the polynomials $P_W(V^0_\lambda,q,y)$(explicitly computed in \cite{GnS}), and then we prove by induction that these closed formulae satisfy Stembridge's recursive relations, specialized in $-q$ and $q^2$. 
In \cite{SDT} we used the minuscule recurrence to prove Reeder's Conjecture in type $B$. A different approach is needed to obtain a proof in type $C$, where Stembridge's quasi minuscule recurrence is more efficient. 
Our choice of using minuscule recurrence in type $D$ comes from a computational reason: in $D_n$ the fundamental weight $\omega_1$ is minuscule and it is consequently easier to reduce the polynomial coefficients to a simpler form.
\section{Small Representations and their Zero Weight Spaces}
We recall that an irreducible representation $V_\lambda$ is \emph{small} if its weight is small.
We refer to \cite{AHJR} and \cite{R4} for a complete classification of small representations and their zero weights spaces as $W$ representations. 

As in the case of hyperoctahedral group $B_n$ examined in \cite{SDT}, the irreducible representations of the Weyl group $D_n=S_n\ltimes \left(\Z/2\Z\right)^{n-1}$
 are encoded by pairs of partitions
$(\nu,\mu),\nu\vdash k,$ $\mu\vdash h, h+k=n$. In the $B_n$ case, the irreducible representations can be realized as $\pi_{\nu,\mu}={\mathrm{ Ind}}_{B_k\times B_h}^{B_n} \pi'_\nu\times \pi''_\mu$,  where, if $\pi_\tau$ is the irreducible $S_p$-module attached to 
$\tau\vdash p$, and $\varepsilon_q$ is the sign representation of  $\left(\mathbb Z_2\right)^q$,  we have $(\pi'_\nu)_{|S_k}=\pi_\nu\,, (\pi'_\nu)_{|\mathbb (Z_2)^k}=1_k, (\pi''_\nu)_{|S_h}=\pi_\mu, 
(\pi''_\mu)_{|\mathbb (Z_2)^h}=\varepsilon_h$. All the irreducible representations of the Weyl group of type $D_n$ can be then obtained restricting to $S_n \ltimes \left(\Z/2\Z\right)^{n-1}$ a representation of the form $\pi_{\nu,\mu}$. If $\nu \neq \mu$, the representation $\widetilde{\pi}_{\nu,\mu}:=\mathrm{Res}^{B_n}_{D_n}\pi_{\nu,\mu}$ remains irreducible, otherwise $\mathrm{Res}^{B_m}_{D_n}\pi_{\nu,\nu}$ splits in two non isomorphic irreducible components that we are going to denote with $\widetilde{\pi}_{\nu,\mu}^I$ and $\widetilde{\pi}_{\nu,\nu}^{II}$.
A complete description of zero weight spaces for small representations in type $D$ is displayed in Table 1.
\begin{table}[h!]\label{smallweightsD}
  \begin{center}
   \caption{Zero weights space of small representations: Type D}
   \begin{tabular}{c|c}
    \textbf{Small Representation} & \textbf{Zero Weight Space } \\
    Highest weight & $(\alpha , \beta)$ description \\
    \hline
    $\omega_{2i}$, $i<\frac{n-1}{2}$ &  $((n-i),(i))$ \\
    $2\omega_{n-1}, \, 2\omega_{n}$ \quad ($n$ \, even) &  $((\frac{n}{2}),(\frac{n}{2}) )_{II}, \;  ((\frac{n}{2}),(\frac{n}{2}) )_{I}$ \\
      $\omega_{n-1}+\omega_{n}$ \quad ($n$ \, odd) &  $((\frac{n+1}{2}),(\frac{n-1}{2}) )$ \\
          $2\omega_1$ &  $((n-1,1), \emptyset)$ \\
    $\omega_1 + \omega_{2i+1}$ , $i<\frac{n-1}{2}$ & $((n-i-1,1),(i)) \oplus ((n-i-1),(i,1))$\\
    $\omega_1 + \omega_{n-1}+\omega_{n}$ \quad ($n$ even) & $((\frac{n}{2},1),(\frac{n-2}{2})) \oplus ((\frac{n}{2}),(\frac{n-2}{2},1))$\\
    $\omega_1 + 2\omega_{n-1}, \, \omega_1 + 2\omega_{n-1} $ \quad ($n$ \, odd)  & $((\frac{n-1}{2},1),(\frac{n-1}{2}))$
\end{tabular}
\end{center}
\end{table}

Let us denote by $\bar{P}_W$ the polynomial $P_W(V^0_\lambda,q,y)$ divided by $\prod_{i=1}^n\left(1-q^{m_i+1}\right)$, where $m_1, \dots, m_n$ are the exponents of Weyl group $W$. We will determine explicit expressions for the polynomial $\bar{P}_W$.

We will encode partitions $\lambda=(\lambda_1, \geq \lambda_2, \dots, \geq, \lambda_n)$ by Young diagrams, displayed in the English way. Here $h(ij)$, $c(ij)$ are the hook length and the content of the box $(ij)$ respectively; $|\lambda|$ and $n(\lambda)$ will denote the quantities $\sum_{i=1}^n \lambda_i$ and $\sum_{i=1}^n(i-1)\lambda_i$. We recall the following results about the $\bar{P}_W$ polynomials for the hyperoctahedral group:
\begin{theorem}[\cite{GnS}, Proposition 3.3]
 Let $\pi_{\alpha, \beta}$ be the irreducible representation of the Weyl group $B_n$ indexed by the pair of partitions $(\alpha, \beta)$.
 \begin{equation}\label{fromabrutatau}
\bar{P}_{B_n}(\pi_{\alpha,\beta}; q,y)= q^{2n(\alpha)+2n(\beta)+|\beta|} \prod_{(i,j) \in \alpha} \frac{1+yq^{2c(ij)+1}}{1-q^{2h(ij)}} \prod_{(i,j) \in \beta} \frac{1+yq^{2c(ij)-1}}{1-q^{2h(ij)}}.
\end{equation}
\end{theorem}
We remark that the reflection representation of the hyperoctahedral group restrict to the reflection representation of the group $S_n \ltimes \left(\Z/2\Z\right)^{n-1}$. Moreover, denoting by $\mathrm{sg(n)}$ the sign representation of $B_n$, the representation $\pi_{\alpha, \beta}$ and $\pi_{\beta, \alpha}= \pi_{\alpha, \beta} \otimes \mathrm{sg(n)}$ restrict to the same irreducible representation $\widetilde{\pi}_{\alpha,\beta}$. As an immediate consequence, it is possible to obtain the following relations between the $P_W$ polynomials in type $B$ and in type $D$:
\[P_D(\widetilde{\pi}_{\alpha,\beta}; q,y)= P_B(\pi_{\alpha,\beta}; q,y) + P_B(\pi_{\beta,\alpha}; q,y)\]
\[P_D(\widetilde{\pi}_{\alpha,\alpha}^I; q,y)=P_D(\widetilde{\pi}_{\alpha,\alpha}^{II}; q,y)= P_B(\pi_{\alpha,\alpha}; q,y)\]
Now we want to rearrange the above formulae in a more useful way. Let $\lambda$ be a partition, we define
\[ S_{\lambda}(q)=(q^n+1)\prod_{(ij) \in \lambda} ( 1+yq^{2c(ij)+1}) \qquad R_{\lambda}(q)=(q^n+1)\prod_{(ij) \in \lambda} ( q+yq^{2c(ij)})\]
In the case of $\lambda=(k,1)$ or $\lambda=(k)$, some nice relations hold: 
\[S_{(k,1)}=S_{(k)}\frac{(q+y)}{q}, \qquad R_{(k,1)}=R_{(k)}\frac{(q^3+y)}{q^2}, \qquad R_{(k)}=q^{k-1}S_{(k-1)}(q+y).\]
Now it is possible to express the $\bar{P}_D$ polynomials for zero weight space representations in Table 1 in a more compact way. Set $H(\lambda)= \prod_{(ij)\in \lambda} \left(1-q^{4h(ij)}\right) $, where $(ij)$ ranges between the boxes of $\lambda$. 
\begin{align*}
 \bar{P}_D(\widetilde{\pi}_{(n-1,1),\emptyset};  q,y)= \frac{q^2\left[R_{(n-1,1)}+S_{(n-1,1)} \right]  }{ H(n-1,1)},
\end{align*}
\begin{align*}
\bar{P}_D(\widetilde{\pi}_{(k),(n-k)}; q,y)&= \frac{S_{(k)}R_{(n-k)}+R_{(k)}S_{(n-k)}     }{H(k) \, H(n-k)\,(1+q^n)}\\&= \frac{q^{k-1}S_{(k-1)} S_{(n-k-1)}(q+y)(q^{n-2k}+1)(1+yq^{n+1})}{H(k) \, H(n-k)\,(1+q^n)},
\end{align*}
\begin{align*}
 \bar{P}_D(\widetilde{\pi}_{(k,1),(n-k-1)} \oplus  &\widetilde{\pi}_{(n-k-1,1),(k)};  q,y)= \\&= \frac{q^2\left[S_{(k)}R_{(n-k-1,1)}+R_{(k)}S_{(n-k-1,1)} \right]  }{H(k) \, H(n-k-1,1)\,(1+q^n)} + \frac{q^2\left[S_{(n-k-1)}R_{(k,1)}+R_{(n-k-1)}S_{(k,1)}\right]    }{H(n-k-1) \, H(k,1)\,(1+q^n)}\\
 &=\frac{q^2 S_{(k-1)}S_{(n-k-2)}(q+y)}{H(k+1)\, H(n-k)\,(1+q^n) \, (1-q^2)} Q(n,k).
\end{align*}
Here we set 
\[P(n,k,q,y)=q^{k-2}(q+y)(1+yq^{2(n-k-2)+1})+q^{n-k-4}(q^3+y)(1+yq^{2(k-1)+1})\]
and
\[Q(n,k,q,y)=P(n,k,q,y)(q^{2(k+1)}-1)(q^{2(n-k-1)}-1) + P(n,n-k-1,q,y)(q^{2k}-1)(q^{2(n-k)}-1).\]

\section{Even Orthogonal Algebras}
We recall now the realization of root system $D_n$ as exposed in \cite{Bou}. In an euclidean vector space with basis $\{e_1, \dots, e_n\}$, consider the set $\Phi_{D_n}$ of vectors $\{\pm e_i \pm e_j\}_{i \neq j, i,j \leq n}$. We fix a positive system of roots choosing the vectors $\{e_i \pm e_j\}_{i\neq j, i,j \leq n}$; according to such a description, the fundamental weights are $\omega_i= e_1+ \dots + e_i$ if $i < n-1$, $\omega_{n-1}= \frac{1}{2}(e_1+ \dots - e_n)$ and $\omega_n= \frac{1}{2}(e_1+ \dots + e_n)$. Moreover the coweight $e_1$ is quasi minuscule and the Weyl vector is $\rho =  \sum (n-j+1) e_j$. 
We have three different families of formulae for $P_W(V_{\lambda}^0, q^2,q)$, depending on the reducibility of representation $V_\lambda^0$.
We will  denote the polynomials $P_D(\widetilde{\pi}_{(k),(n-k)}; q^2,q)$ by $\textbf{C}_{k,n}$, the polynomial $P_D(\widetilde{\pi}_{(n-1,1), \emptyset}; q^2,q)$  by $\textbf{C}_{2|0,n}$ and the polynomials $P_D(\widetilde{\pi}_{(k,1),(n-k-1)} \oplus  \widetilde{\pi}_{(n-k-1,1),(k)};  q^2,q)$ by $\textbf{C}_{2|k,n}$. Similarly we will denote by $C_{k,n}$ and $C_{2|k,n}$ (or simply by $C_{k}$ and $C_{2|k}$ when the context is clear) the rational functions appearing in the Stembridge's recurrences, associated respectively to weight of the form $\omega_{2k}$ and $\omega_1+\omega_{2k+1}$, and specialized in $-q$ and $q^2$.
The following relations between the polynomials $\textbf{C}_{k,n}$ and $\textbf{C}_{2|k,n}$ hold:
\begin{equation}\label{ricorsionePW}
\textbf{C}_{k+1,n}=\textbf{C}_{k,n} \frac{q^2(1+q^{2(n-2k-2)})(1+q^{4k-1})(1-q^{4n-4k})}{(1+q^{2n-4k})(1+q^{4n-4k-5})(1-q^{4k+4})}, 
\end{equation}
\begin{equation}\label{T2kk}
 \textbf{C}_{2|k,n}=
 \textbf{C}_{k,n}\frac{Q(n,k,q^2,q)}{q^{2k-6}(1-q^4)(1-q^{4k+4})(1+q^{2n-1})(1+q^{2n-4k})(1+q^{4n-4k-5})} 
\end{equation}
\begin{rmk}\label{2Paa}
Observe that if we set $n=2k$ in (\ref{ricorsionePW}), the ratio $T_k^n$ between $\textbf{C}_{k,n}$ and $\textbf{C}_{k-1,n}$ becomes equal to 
\[\frac{2\cdot q^2(1+q^{4k-5})(1-q^{4k+4})}{(1+q^{4})(1+q^{4k-1})(1-q^{4k})}\]
This is coherent with the fact that, if $n$ is even, the specialization of $P_D(\widetilde{\pi}_{(k),(n-k)};q^2,q) $ in $k=\frac{n}{2}$ is equal to  $2 P_D(\widetilde{\pi}^I_{(\frac{n}{2}),(\frac{n}{2})};q^2,q)$.
\end{rmk}
\subsection{Weights of the form $\omega_{2k}$, $2 \omega_{n-1}$,  $2 \omega_n$ if $n$ is even and $\omega_{n-1}+\omega_{n}$ if $n$ is odd.}
 We consider the recurrence (\ref{ricorsionemin}) and make the evaluations $q \rightarrow -q$ and $t \rightarrow q^2$, obtaining
\begin{equation}\label{St2}
\sum_{i=1}^lC_{w_i\lambda}\sum_{j=1}^k (q^{-2(\rho,w_i\psi_j)}+q^{1+2(\rho,w_i\psi_j)})=0.
\end{equation}
Writing all the $C_\mu$ in their reduced form the recurrence can be rewritten as 
\begin{equation}\label{redSt2}
\sum_{\mu \leq \lambda}\Lambda_\mu^{\lambda,n}(q)C_\mu(q) = 0,
\end{equation}
for some coefficients $\Lambda_\mu^{\lambda, n}(q)$.
Our purpose is to make more explicit these coefficients. 
\begin{rmk}\label{INVOLUZIONETOTALE}
Set $\Omega_{\mu }^{\lambda, \; n}=\{ w \lambda  \;| w \in W/W_\lambda, \;  w \lambda + \rho \sim \mu+\rho \}$. The sign change on the $n$-th coordinate $\epsilon_n$ induces a bijection between $\Omega_{\mu }^{2 \omega_{n-1}, n}$ and $ \Omega_{\mu }^{2 \omega_{n}, n}$. Observe now that the general coefficient $\Lambda_\mu^{\lambda, n}(q)$ does not depend on the last coordinate of the weights in  $\Omega_{\mu }^{\lambda, \; n}$ because the last coordinate of $\rho$ is equal to 0. As a consequence $\Lambda_\mu^{2\omega_{n-1}, n}(q)=\Lambda_\mu^{2\omega_{n}, n}(q)$ for all $\mu \leq \lambda$ and then $C_{2 \omega_{n-1}}(q)=C_{2 \omega_{n}}(q)$. 
Coherently with the above notation, if $n$ is even we will denote the polynomial $C_{2 \omega_{n}}(q)$ by $C_{\frac{n}{2}}$.
\end{rmk}
 Let us observe that, if $\lambda=\omega_{2k}$ (resp. $\lambda=2\omega_{n}$ and $\lambda=\omega_{n-1}+\omega_{n}$) then the non zero integral dominant weights smaller than $\lambda$ are of the form $\omega_{2i}$ with $i < k$ (resp. $i < \lfloor n/2\rfloor$). 
 For brevity we will denote the coefficient $\Lambda_{\omega_{2i}}^{ \omega_{2k}, n}(q)$ and the set $\Omega_{\omega_{2i}}^{ \omega_{2k}, n}$  (resp. $\Lambda_{\omega_{2i}}^{ 2\omega_{n}, n}(q)$  and $\Omega_{\omega_{2i}}^{ 2\omega_{n}, n}$,  $\Lambda_{\omega_{2i}}^{ \omega_{n-1}+\omega_{n}, n}(q)$  and $\Omega_{\omega_{2i}}^{ \omega_{n-1}+\omega_{n}, n}$ )  by $\Lambda_i^{k,n}$  and $\Omega_i^{k,n}$ (resp. by $\Lambda_i^{\frac{n}{2},n}$ and $\Omega_i^{\frac{n}{2},n}$, $\Lambda_i^{\frac{n-1}{2},n}$ and $\Omega_i^{\frac{n-1}{2},n}$). 
It is immediate to show that the only contributions to the coefficient $\Lambda_k^{k, n}$ come from the case $w_i=id$  and consequently $\Lambda_k^{k, \,n}$ must be equal to
\begin{equation}\label{coefficientek}
\Lambda_k^{k, \,n}=\sum_{i=1}^{2k}\frac{(1+q^{4(n-i)+1})}{q^{2(n-i)}}=\frac{(q^{4k}-1)(q^{4(n-k)-1}+1)}{q^{2(n-1)}(q^2-1)}.
\end{equation}
It is more difficult to obtain closed formulae for the generic $\Lambda_{h}^{k,n}$, however some nice recurrences hold.
\begin{rmk}\label{omegah+zero}
A direct inspection shows that the weights giving non zero contribution to $\Lambda_{h}^{k,n}, k>h>0$ are of the form $e_1+\dots+e_{2h}+\mu$, where $\mu$ has the first $2h$ coordinates equal to 0. Considering the immersion of $D_{n-2h} \rightarrow D_n$ induced by Dynkin diagrams, this means that $\mu$ can be contracted to a weight in $\Omega_{0}^{k-h,n-2h}$. By abuse of notation, we denote this contraction process writing~$\mu\in\Omega_{0}^{k-h,n-2h}$.
\end{rmk}
As a consequence of the above Remark, the following relation holds \footnote{Here we set by convention $\Omega_{0}^{1,2}=\{(-1,-1)\}$, $\Omega_{0}^{1,3}=\{(-1,1,0),(0,-1,-1),(0,-1,1\}$ and coherently $\Lambda_0^{1,2}=r(2,q)$, $\Lambda_0^{1,3}=r(3,q)+2r(2,q)$(c.f.r. Equation (\ref{L002k}))  }:
\begin{equation}\label{Lh}
\Lambda_h^{k, \, n}= (-1)^{k-h}\Lambda_h^{h,n} |\Omega_0^{k-h, n-2h}| + \Lambda_0^{k-h, \, n-2h}.  \end{equation}
We reduced to computing of the coefficients $\Lambda_0^{k, \, n}$.
For the weights conjugated to 0 some results similar to the ones proved in \cite{SDT} for $B_n$ and $C_n$ hold:
\begin{lemma}\label{zerocontribution1}
Let $n$ be even. Set $\lambda=2 \omega_n$ and let $w\in W$ be such that $w \lambda$ is conjugated to $0$. 
\begin{enumerate}
 \item $w \lambda$ is equal to $(-1, 1, \dots , -1,1)$ if $n/2$ is even and to $(-1, 1, \dots , -1,-1)$ if $n/2$ is odd.
 \item There exists an element $\sigma \in W$ of sign $\mathrm{sg}(\sigma)=n/2$ such that $\sigma (w \lambda + \rho)= \rho$.
\end{enumerate}
\end{lemma}
\begin{lemma}\label{zerocontribution2}
Set $\lambda=\omega_{2k}$, $2k<n$ or $\lambda=\omega_{n-1}+\omega_{n}$, $n=2k+1$ and let $w\in W$ be such that $w \lambda$ is conjugated to $0$. Then $w \lambda$ is of one of the following form:
\begin{enumerate}
 \item The $2k$ non zero coordinates of $w \lambda$ are pair of consecutive coordinates $((w\lambda)_{(j)},(w\lambda)_{(j)+1})$  of the form $(-1,1)$.
 \item There are $2(k-1)$ non zero coordinates that are pair of consecutive coordinates $((w\lambda)_{(j)},(w\lambda)_{(j)+1})$  of the form $(-1,1)$ and the latter two are equal to $-1$.
\end{enumerate}
In both cases there exists an element $\sigma \in W$ of sign $\mathrm{sg}(\sigma)=k$ such that $\sigma (w \lambda + \rho)= \rho$.
\end{lemma}
This produce an explicit formula for the number of weights $\Omega_0^{\lambda, n}$:
\begin{equation*}\label{CardG0}
 \displaystyle{
| \Omega_0^{\lambda, n}| =
\left\{
	\begin{array}{lll}
	1 & \mbox{if } \lambda=2\omega_n \mbox{ and }   n \,  \mbox{ is even or } \lambda=0 \\
    \frac{n}{k}\binom{n-k-1}{k-1}& \mbox{if } \lambda=\omega_{2k}  \mbox{ and }   2k < n \mbox{ or } \lambda=\omega_{n-1}+\omega_{n}  \mbox{ and } n=2k+1.\\
	\end{array}
\right.}
\end{equation*}
Furthermore, the coefficient $\Lambda_0^{\frac{n}{2},n}$ can be easily computed. Set 
\[r(n,q)=q^{2(n-1)}-q^{-2(n-1)+1}+q^{-2(n-2)}+q^{2(n-2)+1}=(q+1)\left(q^{2n-3}+q^{-2n+3}\right)\]
then 
\begin{equation}\label{L002k}
 \Lambda_0^{\frac{n}{2},n} = (-1)^{\frac{n}{2}}\sum_{i=1}^{\frac{n}{2}} r(2i, q)=(-1)^{\frac{n}{2}} r(n,q) - \Lambda_0^{\frac{n-2}{2}, n-2}.
\end{equation}
In the general case, producing explicit formulae for $\Lambda_0^{k,n} $ is more complicated.
As a consequence of Lemma \ref{zerocontribution2} we obtain a case by case analysis of the weights in $\Omega_{0}^{k, \, n}$ that leads us to some recursive expressions for  $\Lambda_0^{k,n}$. Let $\mu$ be a weight in $\Omega_{0}^{k, \, n}$:
\begin{itemize}
\item[\emph{Case 1:If $n=2k+1$}] then $\mu$ must be of the form $-e_1+e_2+\nu$ with $\nu \in \Omega_0^{k-1,2k-1}$ or $\mu=(0, \mu_2, \dots, \mu_n)$ where $\mu'=(\mu_2, \dots, \mu_n) \in \Omega_0^{k,2k}$ or $\epsilon_n \mu' \in \Omega_0^{k,2k}$,
\item[\emph{Case 2: If $n\neq 2k, 2k+1$}] then $\mu$ must be of the form $-e_1+e_2+\nu$ with $\nu \in \Omega_0^{k-1,n-2}$ or $\mu=(0, \mu_2, \dots, \mu_n)$ where $\mu'=(\mu_2, \dots, \mu_n) \in \Omega_0^{k,n-1}$. 
\end{itemize}
We obtain the recursive relations:
 \begin{equation}\label{L02k1}\Lambda_{0}^{k, \, 2k+1}= (-1)^k r(2k+1,q)|\Omega_0^{k-1,2k-1}|- \Lambda_0^{k-1,2k-1} + 2 \Lambda_0^{k,2k}, \end{equation}
  \begin{equation}\label{L0n}\Lambda_{0}^{k, \, n}= (-1)^k r(n,q)|\Omega_0^{k-1,n-2}|- \Lambda_0^{k-1,n-2} + \Lambda_0^{k,n-1}. \end{equation}
The above recurrences between the coefficients allow us to reduce the triangular system given by Stembridge's relations as described in the following proposition, that we prove in Section 5.
Set 
 \begin{equation*}\label{bkn}
 \displaystyle{
b_{k,n} =
\left\{
	\begin{array}{lll}

	\frac{(q^{4k}-1)}{q^{2k-1}(q-1)}& \mbox{if }    n=2k, \\
 	& \\
  \frac{(q^{2n}-1)(q^{2(n-2k)}+1)}{q^{2(n-k)-1}(q-1)}& \mbox{otherwise. } \\


	\end{array}
\right.}
\end{equation*}
\begin{proposition}\label{redrecursion}
 Let $R_i$ be the recurrence for $C_{i}$ written in reduced form. Then there exist a family of integers $\{A_i^{k, n}\}_{i \leq k}$ such that 
 \begin{equation}\label{Riduzionebella}\sum_{i=1}^k A_i^{k,n}R_i = \Lambda_k^{k, n} C_{k}- \sum_{i=1}^k b_{i,n-2(k-i)}C_{k-i}\end{equation}
\end{proposition}
Now, we apply an inductive reasoning. We recall that the base case of weight $\omega_2$ can be obtained comparing $P_D(V_{\omega_2}^0, q^2, q)$ with the formulae proved by Stembridge in \cite{Stembridge}.
Observe now that the following relations between the $b_{k,n}$ hold:
\begin{equation*}b_{i,n-2(k-i)}-b_{i-1,n-2(k-i)}=-\frac{(q+1)(q^{2(n-2k+1)}-1)}{q^{2(n-2k)+1}} \cdot \frac{(q^{2(n-2(k-i))}-1)}{q^{2i}},\end{equation*}
\begin{equation*}2 b_{i,2i}-b_{i-1,2i}=-\frac{(q+1)(q^{2}-1)}{q} \cdot \frac{(q^{4i}-1)}{q^{2i}}.\end{equation*}
If $k<\frac{n}{2}$, considering the difference between the recurrences of the form (\ref{Riduzionebella}) for $C_k$ and $C_{k-1}$ we obtain:
\begin{align*}
 \Lambda_k^{k,n}C_k - \left(b_{1,n-2(k-1)}+\Lambda_{k-1}^{k-1,n}\right)C_{k-1} &=   \sum_{i=2}^k \left[b_{i,n-2(k-i)}-b_{i-1,n-2(k-i)}\right]C_{k-i} \\
 &=-\frac{(q+1)(q^{2(n-2k+1)}-1)}{q^{2(n-k)+1}} \sum_{i=0}^{k-2}q^{2i}(q^{2(n-2i)}-1)C_i
\end{align*}
Set
\[C(k,n,q)=\Lambda_k^{k,n}T_k^n - b_{1,n-2(k-1)}-\Lambda_{k-1}^{k-1,n},\]
\[D(k,n,q)=-\frac{(q+1)(q^{2(n-2k+1)}-1)}{q^{2(n-k)+1}}.\]
The conjecture in this case is then reduced to prove that 
\begin{equation}\label{polidCk}
 \frac{C(k,n,q)}{D(k,n,q)} C_{k-1}=\left[q^{2(k-2)} \left(q^{2(n-2(k-2))}-1\right)+\frac{C(k-1,n,q)}{D(k-1,n,q)}\right]C_{k-2}.
\end{equation}
and then to check a polynomial identity. This can be easily computed by a symbolic algebra software. We checked the above equality using SageMath \cite{Sage}. 
In the case of $2 \omega_n$ with $n=2k$, we have to consider 2 times (\ref{Riduzionebella}) for $C_\frac{n}{2}$ minus the recurrence for $C_{k-1}$, obtaining 
\begin{align*}
 \Lambda_k^{k,n}2 C_n - \left(2 b_{1,2}+\Lambda_{k-1}^{k-1,n}\right)C_{k-1} &=   \sum_{i=2}^k \left[2 b_{i,2i}-b_{i-1,2i}\right]C_{k-i} \\
 &=-\frac{(q+1)(q^{2}-1)}{q^{n+1}} \sum_{i=0}^{k-2}q^{2i}(q^{2(n-2i)}-1)C_i,
\end{align*}
and the computation is exactly the same as in the previous case. Observe that our choice of take 2 times recurrence (\ref{Riduzionebella}) for $C_\frac{n}{2}$ in the above computation is coherent with the Remark \ref{2Paa}.
\subsection{Weights of the form $\lambda=\omega_1 + \omega_{2k+1}$, $\lambda=\omega_1 + \omega_{n-1}+\omega_n$ if $n$ is even and $\lambda=\omega_1 + 2\omega_{n-1}, \omega_1 + 2\omega_n$ if $n$ is odd } We remark that in this case $W_\lambda \cdot e_1=\{e_1\}$ and then the minuscule recurrence becomes more explicit.
 We consider the recurrence (\ref{ricorsionemin}) and make again the evaluations $q \rightarrow -q$ and $t \rightarrow q^2$, obtaining
\begin{equation}\label{St3}
\sum_{i=1}^lC_{w_i\lambda} \frac{( 1-q^{1+4(\rho,w_ie_1)})}{q^{2(\rho,w_ie_1)}}=0.
\end{equation}
The recurrence (\ref{St3}) can be written in reduced form as:
\begin{equation*}
  R_k  \; :  \; \sum_{i=0}^k C_{2|i}\Lambda^{2|k,n}_{2|i}+\sum_{i=0}^{k+1/k} C_{i}\Lambda_i^{2|k, \, n}=0,
 \end{equation*}
 for some coefficients $\Lambda^{2|k,n}_{2|i}$ and $\Lambda^{2|k,n}_{i}$ that we are going now to analyze more closely. We underline that the index of the second sum goes from 0 to $k$ if $n=2k+1$ and to $k+1$ if $n>2k+1$.
In the special case $n=2k+2$ the recurrence can be displayed as follows:
\begin{equation*}
  R_k  \; :  \; \sum_{i=0}^k C_{2|i}\Lambda^{2|k,n}_{2|i}+\sum_{i=0}^{k} C_{i}\Lambda_i^{2|k, \, n}+ C_{k+1}\left[\Lambda_{2 \omega_n}^{2|k, \, n} +  \Lambda_{2 \omega_{n-1}}^{2|k, \, n}\right]=0.
 \end{equation*}
Using an argument similar to Remark \ref{INVOLUZIONETOTALE} we have $\Lambda_{2 \omega_{n}}^{2|k, \, n} =  \Lambda_{2 \omega_{n-1}}^{2|k, \, n}$ and then the recurrence becomes: 
\begin{equation*}
  R_k  \; :  \; \sum_{i=0}^k C_{2|i}\Lambda^{2|k,n}_{2|i}+\sum_{i=0}^{k} C_{i}\Lambda_i^{2|k, \, n}+ 2C_{n}\Lambda_{k+1}^{2|k, \, n} =0.
 \end{equation*}
 Moreover, by the same reasoning of Remark \ref{INVOLUZIONETOTALE}, if $n=2k+1$ we obtain that $C_{\omega_1 + 2\omega_{n-1}}(q)=C_{\omega_1 + 2\omega_{n}}(q)$ and then we reduce to compute only $C_{\omega_1 + 2\omega_{n}}(q)$ (that for brevity we will denote with $C_{2|k, \, 2k+1}$ ).
As in the previous case we want to determine closed formulae and recursive relations for the coefficients that allow us to reduce the system. 
\begin{rmk}\label{rho+peso1k}
Similarly to what we observed in the case $\omega_{2k}$, the only contributions to $\Lambda^{2|k, n}_{2|h}$ come from weights of the form $2e_1+e_2 + \dots + e_{2h+1} + \mu$ with $\mu \in \Omega_{0}^{k-h,n-2h-1}$.
\end{rmk}
Set \[s(n,q)=\frac{(1-q^{4n+2})}{q^{2n}},\]
by previous remark we obtain:
\begin{equation}\label{coeff2h}
\Lambda^{2|k, n}_{2|h}=(-1)^{k-h}s(n-1,q)\left|\Omega_{0}^{k-h,n-2h-1}\right|.
\end{equation}
Coherently with the notation of the previous section, we will denote by $\Omega_{2|h}^{2|k,n}$ (resp. $\Omega_{h}^{2|k,n}$)  the set $\Omega_{\omega_1+\omega_{2h+1}}^{\omega_1+\omega_{2k+1},n}$ (resp. $\Omega_{\omega_{2h}}^{\omega_1+\omega_{2k+1},n}$). It is immediate to show that the weights in $\Omega_{k+1}^{2|k,n}$ are all of the form $e_1 + \dots + e_{j-2} + 2e_{j} + \dots + e_{2(k+1)} $ for $j\leq 2(k+1)$.
As a consequence, we obtain the following formula:
\begin{equation}\label{Lkk+1}
 \Lambda^{2|k, n}_{k+1} = - \sum_{j=2}^{2k+2}s(n-j,q)=-\frac{(q^{4(n-k-1)-2}-1)(q^{4(k+1)-2}-1)}{q^{2(n-2)}(q^2-1)}.
\end{equation}
\begin{rmk}{\textbf{Reeder's Conjecture for $2\omega_1$.}}
Formulae (\ref{coeff2h}) and (\ref{Lkk+1}) describe explicitly the coefficients $\Lambda_{2|0}^{2|0,n}$ and $\Lambda_1^{2|0,n}$. By direct inspection it is possible to prove that $\Omega_{0}^{2|0,n}=\{-2e_{n-1}\}$ and consequently $\Lambda_{0}^{2|0,n}=s(-1,q)=-s(0,q)$. Now $R_0$ is explicit and the Reeder's Conjecture in this case can be proved comparing the obtained expression for $C_{2|0,n}$ with the formula for  $P_D(\widetilde{\pi}_{(n-1,1),\emptyset};  q^2,q)$.\footnote{Observe that, to prove inductively the Conjecture, it is also needed to check explicitly the cases $\omega_1+\omega_3$ in $D_4$ and $\omega_1+\omega_5$ in $D_5$. In these cases to list directly the weights involved in the computation of coefficients is quite simple. It is possible to achieve a description coherent with the general analysis in the paper setting by convention $\Omega_{0}^{2|0,n}=\{ (-2,0)\}$ for \emph{every} $n>1$ and $\Omega_{0}^{2|1,3}=\{(-1,-1,2), (-2,1,-1)\}$, $\Omega_{1}^{2|1,3}=\{(1,-1,2)\}$. Coherently $\Lambda^{2|1,3}_{0}=s(0,q)+s(-2,q)$ and $\Lambda^{2|k-1, 2(k-1)+1}_{h-1}=-s(0,q)$. }
\end{rmk}
To find a recursive expansion of the coefficients $\Lambda^{2|k, n}_h$ for $k\geq 1$ we have to consider the cases $n=2k+1,2k+2$ and $n$ generic in three different ways. Firstly suppose $k \geq h>1$.
If $n=2k+1$ a weight $\mu \in \Omega_h^{2|k,n}$  must be of the form $\omega_{2}+\mu'$ with $\mu' \in \Omega_{h-1}^{2|k-1,n-2}$, obtaining 
\[\Lambda^{2|k, n}_h=\Lambda^{2|k-1, 2(k-1)+1}_{h-1}.\]
The case $n=2k+2$ needs a closer analysis. A weight $\mu \in \Omega_h^{2|k,n}$ can be of one of the following forms:
\begin{itemize}
 \item[Case 1] $\mu = e_1+e_2 + \mu'$ with $\mu' \in \Omega_{h-1}^{2|k-1,n-2}$,
 \item[Case 2] $\mu = 2e_2 + \mu'$ with $\mu' \in \Omega_{h-1}^{k,n-2} \cup \epsilon_n \Omega_{h-1}^{k,n-2}$ or equivalently $\mu_1=0$ and $\tilde{\mu}=(\mu_2, \dots, \mu_n) \in \Omega_{2|h-1}^{2|k,n-1} \cup \epsilon_n\Omega_{2|h-1}^{2|k,n-1}$,
 \item[Case 3] $\mu=(1,0,2,1,\mu_5, \dots, \mu_n)$ where $\tilde{\mu}=(\mu_5, \dots, \mu_n) \in \Omega_{h-2}^{k-1,2k-2} \cup \epsilon_n\Omega_{h-2}^{k-1,2k-2}$.
\end{itemize}
This analysis can be translated in the following relation:
\[\Lambda^{2|k, 2k+2}_h=\Lambda^{2|k-1, 2k}_{h-1}-2\Lambda^{2|k, 2k+1}_{2|h-1}-2(-1)^{k-h+1}s(2k-1,q).\]
The case of $n$ generic is completely analogous to the previous one, except for the fact that $\Omega_{h-1}^{k,n-2} \cup \epsilon_n \Omega_{h-1}^{k,n-2}=\Omega_{h-1}^{k,n-2}$ and 
$\Omega_{h-2}^{k-1,n-4} \cup \epsilon_n\Omega_{h-2}^{k-1,n-4}=\Omega_{h-2}^{k-1,n-4}$
leading to the recurrence:
\[\Lambda^{2|k, n}_h=\Lambda^{2|k-1, n-2}_{h-1}-\Lambda^{2|k, n-1}_{2|h-1}-(-1)^{k-h+1}s(n-3,q)\left|\Omega_{0}^{k-h+1, \, n-2h}\right|.\]
The case $h=1$ must be considered differently. If $n=2k+1$, the elements in $\Omega_{1}^{k,2k+1}$ are of the form:
\begin{itemize}
 \item[Case 1] $\mu = \omega_2 + \mu'$ with $\mu' \in \Omega_0^{2|k-1, \, n-2}$,
 \item[Case 2] $\mu = (1,-1,2, \mu_4, \dots, \mu_n)$ with $\tilde{\mu}=(\mu_4, \dots, \mu_n) \in \Omega_0^{k-1,n-3}$,
\end{itemize}
obtaining the recurrence:
\[\Lambda_1^{2|k,2k+1}=\Lambda_0^{2|k-1,2k-1}-(-1)^{k-1}s(2k-2,q).\]
In the other cases we must consider a third family of weights, i.e. of the form $(0, \mu_2, \dots, \mu_n)$ where $\tilde{\mu}=(\mu_2, \dots, \mu_n)$ is contained in $\Omega_{2|0}^{2|k,2k+1} \cup \epsilon_n \Omega_{2|0}^{2|k,2k+1}$ if $n=2k+2$ and in $\Omega_{2|0}^{2|k,n-1}$ for $n$ generic.
Consequently we obtain:
\[\Lambda_1^{2|k,2k+2}=\Lambda_0^{2|k-1,2k} -2 \Lambda_{2|0}^{2|k, 2k+1}-(-1)^{k-1}s(2k-1,q)\left|\Omega_0^{k-1,2k-1}\right|,\]
\[\Lambda_1^{2|k,n}=\Lambda_0^{2|k-1,n-2} - \Lambda_{2|0}^{2|k, n-1}-(-1)^{k-1}s(n-3,q)\left|\Omega_0^{k-1,n-3}\right|.\]
Finally, we will find recurrences for coefficient $\Lambda_0^{2|k,n}$.
Again, let us start from $n=2k+1$. The only relevant cases are: 
\begin{itemize}
 \item[Case 1] $\mu = (-1,1, \mu_3, \dots, \mu_n)$ with $\tilde{\mu}=(\mu_3, \dots, \mu_n) \in \Omega_0^{2|k-1,n-2}$,
 \item[Case 2] $\mu = (-1,-1,2, \mu_4, \dots, \mu_n)$ with $\tilde{\mu}=(\mu_4, \dots, \mu_n) \in \Omega_0^{k-1,n-3}$,
  \item[Case 3] $\mu = (-2,1,1, \mu_4, \dots, \mu_n)$ with $\tilde{\mu}=(\mu_4, \dots, \mu_n) \in \Omega_0^{k-1,n-3}$,
\end{itemize}
and consequently we obtain: 
\[\Lambda_0^{2|k,2k+1}=-\Lambda_0^{2|k-1,2k-1}+(-1)^{k-1}\left[s(2k-2,q)+s(-2k,q)\right].\]
As in the case $h=1$, if $n\neq 2k+1$ we have to add to the previous list the weights of the form $(0,\mu_2,\dots, \mu_n)$ with $\tilde{\mu}=(\mu_2, \dots, \mu_n)$ contained respectively in $\Omega_{0}^{2|k,2k+1} \cup \epsilon_n \Omega_{0}^{2|k,2k+1}$ if $n=2k+2$ and in $\Omega_{0}^{2|k,n-1}$ for $n$ generic. This yields to the recurrences: 
\[\Lambda_0^{2|k,2k+2}=2 \Lambda_0^{2|k,2k+1}-\Lambda_0^{2|k-1,2k}+(-1)^{k-1}\left[s(2k-1,q)+s(-2k-1,q)\right]\left|\Omega_0^{k-1,2k-1}\right|,\]
\[\Lambda_0^{2|k,n}= \Lambda_0^{2|k,n-1}-\Lambda_0^{2|k-1,n-2}+(-1)^{k-1}\left[s(n-3,q)+s(-n+1,q)\right]\left|\Omega_0^{k-1,n-3}\right|.\]
Set now
\[d_{k,n}:= s(k,q)+s(n-k-1,q)=\frac{q^{2k}(1-q^{2n})(1+q^{2(n-2k-1)})}{q^{2(n-1)}},\]
imposing that the coefficients $\Lambda^{2|k,n}_{2|h}$ must cancel, we are again able to find a nice reduction to the triangular system, proved in Section 5.
\begin{proposition}\label{magic}
Let $R_i$ be the reduced recurrence for $C_{2|i}$. Then there exist 
 a family of integers $\{E_i^{k, n}\}_{i \leq k}$, with $E_k^{k,n}=1$ and $\sum_{i=h}^k(-1)^i \binom{n-i-1}{i}E_{i}^{k n }=0$, such that 
 \begin{equation}\label{n=2k+1}
 \sum_{i=0}^k E_i^{k,2k+1}R_i = \Lambda_{2|k}^{2|k, 2k+1} C_{2|k,2k+1}+ \Gamma_{k}^{2|k, 2k+1} C_{k,2k+1}  - \sum_{j=0}^{k-1}s(k-j,q)C_{j,2k+1},
\end{equation}
 and
 \begin{equation}\label{ngenerico}
\sum_{i=0}^k E_i^{k,n}R_i = \Lambda_{2|k}^{2|k, n} C_{2|k,n}+ \Gamma_{k+1}^{2|k, n} C_{k+1,n} + \Gamma_{k}^{2|k, n} C_{k,n}  - \sum_{j=0}^{k-1}d_{k-j,n-2j}C_{j,n}.
 \end{equation}
 for some suitable coefficients $\Gamma_{k}^{2|k, 2k+1}$, $\Gamma_{k+1}^{2|k, n}$ and $\Gamma_{k}^{2|k, n}$.
\end{proposition}
We remark that it is possible to find a recursive expression similar to (\ref{Lkk+1}) for the coefficient $\Gamma_k^{2|k,n}$ (c.f.r. Section 5, Proposition \ref{magic5S}, computations for $\Gamma_h^{2|k,n}$):
\[\Gamma_1^{2|1,n}=-s(0,q)-s(n-2,q)-s(n-3,q)\]
\begin{align*}
 \Gamma_k^{2|k,n}
 =\Gamma_{k-1}^{2|k-1,n-2}-\left[s(n-2,q)+s(n-3,q)\right]
\end{align*}
A direct inspection shows that $\Gamma_2^{2|2,5}=-\sum_{j=0}^3s(j,q)$ and, starting from the formula for $\Gamma_1^{2|1,n}$, it is possible to obtain the following expressions:
\[\Gamma_k^{2|k,2k+1}=-\sum_{j=0}^{2k-1}s(j,n)=-\frac{(q^{4k}-1)^2}{q^{4k-2}(q^2-1)},  \qquad \Gamma_k^{2|k,n}=-s(0,q)+s(n-2k-2,q)+ \Gamma_{k+1}^{2|k,n}.\]
\subsection{Reeder's Conjecture for $\lambda=\omega_1+\omega_{2k+1}$, $\lambda=\omega_1+\omega_{n-1}+\omega_n$ if $n$ is even and $\lambda=\omega_1+2\omega_{n-1}, \omega_1+2\omega_n$ if $n$ is odd.}
We recall that in Section 4.1 we shown the identity:
\begin{equation}\label{ricorsioneperinduzione}
\frac{C(k,n,q)}{D(k,n,q)}C_{k-1} = \sum_{i=0}^{k-2}q^{2i}(q^{2(n-2i)}-1)C_i.
\end{equation}
Substituting (\ref{ricorsioneperinduzione}) in (\ref{n=2k+1}) we obtain 
\begin{align*}
 \Lambda_{2|k}^{2|k, 2k+1} C_{2|k,2k+1}&=- \Gamma_{k}^{2|k, 2k+1} C_{k,2k+1}  + \sum_{j=0}^{k-1}s(k-j,q)C_{j,2k+1}\\
 &=- \Gamma_{k}^{2|k, 2k+1} C_{k,2k+1}  - \sum_{j=0}^{k-1}\frac{q^{2j}(q^{2(2k+1-2j)}-1)}{q^{2k}}C_{j,2k+1}\\
 &=- \Gamma_{k}^{2|k, 2k+1} C_{k,2k+1}  - \frac{C_{k-1,2k+1}}{q^{2k}} \left[q^{2(k-1)}(q^{6}-1)+\frac{C(k,n,q)}{D(k,n,q)}\right]\\
 &=- \frac{C_{k-1,2k+1}}{q^{2k}} \left[q^{2k}\Gamma_{k}^{2|k, 2k+1} T_k^{2k+1}  + q^{2(k-1)}(q^{6}-1)+\frac{C(k,n,q)}{D(k,n,q)}\right],
 \end{align*}
and analogously, substituting (\ref{ricorsioneperinduzione}) in (\ref{ngenerico})
\begin{align*}
\Lambda_{2|k}^{2|k, n} C_{2|k,n}&=- \Gamma_{k+1}^{2|k, n} C_{k+1,n}  -\Gamma_{k}^{2|k, n} C_{k,n}+\sum_{j=0}^{k-1}d_{k-j,n-2j}C_{j,n}\\
 &=- \Gamma_{k+1}^{2|k, n} C_{k+1,n}  -\Gamma_{k}^{2|k, n} C_{k,n}+\frac{(q^{2(n-2k-1)}+1)}{q^{2(n-k-1)}}\sum_{j=0}^{k-1}q^{2j}(q^{2(n-2j)}-1)C_{j,n}\\
 &=- C_{k,n}\left[\Gamma_{k+1}^{2|k, n} T_{k+1}^n  + \Gamma_{k}^{2|k, n}+\frac{(q^{2(n-2k-1)}+1)}{q^{2(n-k-1)}}\frac{C(k,n,q)}{D(k,n,q)}\right].
\end{align*}
Again proving Reeder's Conjecture is equivalent to a polynomial identity that we checked with SageMath\cite{Sage}.
\section{Proof of Propositions 4.5 and 4.7}
In this Section we are going to denote with $\Theta(n.k,h)$ the quantity $\left|\Omega_h^{k,n}\right|=\left|\Omega_0^{k-h,n-2h}\right|$. 
\begin{rmk}\label{tetaprop}
 \[\Theta(n,k,h)=\Theta(n,k,h+1) + \Theta(n-1,k,h),\]
  \[\Theta(2k+1,k,h)=\Theta(2k+1,k,h+1) + 2\Theta(2k,k,h).\]
\end{rmk}
Set now
\begin{equation}\label{blabla}
A_h^{k, n} =
\left\{
	\begin{array}{lll}
        0  & \mbox{if } h>k \mbox{ or } h\leq 0,\\
		1  & \mbox{if } h=k, \\
		-\sum_{i=h+1}^k(-1)^{i-h} \Theta(n,i,h) A_{i}^{k,n} & \mbox{otherwise.}
	\end{array}
\right.
\end{equation}
\begin{proposition}
 Let $R_i$ be the recurrence for $C_{i}$ written in reduced form. Then
 \begin{equation}\label{RiduzionebellaS5}\sum_{i=1}^k A_i^{k,n}R_i = \Lambda_k^{k, n} C_{k}- \sum_{i=0}^{k-1} b_{k-i,n-2i}C_{i}.\end{equation}
 \proof
 We will use $\Gamma_{h}^{k,n}$ to denote the coefficient of $C_{h}$ in $\sum_{i=1}^k A_i^{k,n}R_i $. 
 Using remark \ref{tetaprop} and definition (\ref{blabla}) it is possible to prove that the integers of the form $A_h^{k,n}$ satisfy some nice iterative properties:
\begin{lemma}\label{propAhkn}
 \begin{enumerate}
  \item $A_h^{k,n}=A_{h-1}^{k-1,n-2}$ $\mathrm{for}$ $h>1$,
  \item $A_{h}^{k,n}=A_{h}^{k,n-1}+A_{h}^{k-1,n-1}$ $\mathrm{if}$ $n \neq 2k, 2k+1$,
  \item $A_{h}^{k,2k}=A_{h}^{k-1,2k-1}$,
  \item $A_{h}^{k,2k+1}=2A_{h}^{k,2k}+A_{h}^{k-1,2k}$.
 \end{enumerate}
\end{lemma}
Using $(1)$ it is easy to prove that $\Gamma_h^{k,n}=\Gamma_0^{k-h,n-2h}$ if $k>h>0$:
\begin{align*}
 \Gamma_h^{k,n} &= \sum_{j=h}^kA_j^{k,n}\Lambda_h^{j,n}\\
 &=\left[\sum_{j=h}^k (-1)^{j-h}\Theta(n,j,h) \, A_j^{k,n}\right] \Lambda_h^{h,n} + \sum_{j=h+1}^k  A_j^{k,n} \Lambda_0^{j-h,n-2h}\\
 &=\sum_{t=1}^{k-h}  A_t^{k-h,n-2h} \Lambda_0^{t,n-2h}=\Gamma_0^{k-h,n-2h}.
\end{align*}
To recover the coefficients of (\ref{RiduzionebellaS5}) we have to compute $\Gamma_0^{k,n}$. In particular we want to prove by inductive reasoning that 
\[\Gamma_0^{k,2k}=-\sum_{j=2}^{k+1}r(j,q), \qquad \Gamma_0^{k,2k+1}=2 \Gamma_0^{k,2k}-r(k+2,q), \qquad \Gamma_0^{k,n}= \Gamma_0^{k,n-1}-r(n-k+1,q).\]
The formula for $\Gamma_0^{2,4}$ can be easily checked by direct computations. We have to distinguish three different cases. 
If $n=2k$ and using $(1)$ and $(3)$ of Lemma \ref{propAhkn} we have:
\begin{align*}
  \Gamma_0^{k,2k}&=\sum_{j=1}^{k} \left[\Lambda_0^{j,n-1}-\Lambda_0^{j-1,n-2}\right]A_j^{k.n} +\sum_{j=1}^k\left[(-1)^{j}\Theta(n,j,1)A_j^{k,n} \right]  r(2k,q)\\
  &=\sum_{j=1}^{k-1} A_j^{k-1,n-1} \Lambda_0^{j,n-1} - \sum_{t=1}^{k-1} A_t^{k-1,n-2} \Lambda_0^{t,n-2}\\
  &=\Gamma_0^{k-1,2(k-1)+1}-\Gamma_0^{k-1,2(k-1)}=^{Ind}\Gamma_0^{k-1,2(k-1)}-r(k+1,q).
 \end{align*}
The computation in the case $n=2k+1$ needs the additional use of $(4)$.  
 \begin{align*}
  \Gamma_0^{k,2k+1}&=\sum_{j=1}^{k-1} \left[\Lambda_0^{j,n-1}-\Lambda_0^{j-1,n-2}\right]A_j^{k,n} +\sum_{j=1}^k\left[(-1)^{j}\Theta(n,j,1)A_j^{k,n} \right]  r(2k+1,q)+2 \Lambda_0^{k,2k}-\Lambda_0^{k-1,2k-1}\\
  &=\sum_{j=1}^{k-1} A_j^{k,n} \Lambda_0^{j,2k}+ 2 \Lambda_0^{k,2k}- \sum_{t=1}^{k-1} A_t^{k-1,n-2} \Lambda_0^{t,n-2}\\
  &=\sum_{j=1}^{k-1} \left[A_j^{k-1,2k}+2A_j^{k,2k}\right] \Lambda_0^{j,2k}+ 2 \Lambda_0^{k,2k}-\Gamma_0^{k-1,n-2}\\
  &=\Gamma_0^{k-1,2k}+2\Gamma_0^{k,2k}-\Gamma_0^{k-1,n-2}\\
  &=^{Ind}2\Gamma_0^{k,2k)}-r(k+2,q).
 \end{align*}
The case of $n$ generic is completely analogous and uses $(2)$ in lemma \ref{propAhkn} instead of $(3)$.
 \begin{align*}
  \Gamma_0^{k,n}&=\sum_{j=1}^{k} \left[\Lambda_0^{j,n-1}-\Lambda_0^{j-1,n-2}\right]A_j^{k,n} +\sum_{j=1}^k\left[(-1)^{j}\Theta(n,j,1)A_j^{k,n} \right]  r(n,q)\\
  &=\sum_{j=1}^{k} \left[A_j^{k-1,n-1}+A_j^{k,n-1}\right] \Lambda_0^{j,n-1}-\Gamma_0^{k-1,n-2}\\
  &=\Gamma_0^{k-1,n-1}+\Gamma_0^{k,n-1}-\Gamma_0^{k-1,n-2}\\
  &=^{Ind}\Gamma_0^{k,n-1}-r(n-k+1,q).
 \end{align*}
 Now it is straightforward to show that $\Gamma_0^{k,n}=b_{k,n}$ and then $\Gamma_h^{k,n}=\Gamma_0^{k-h,n-2k}=b_{k-h,n-2h}$. 
 \endproof
\end{proposition}
Set now $\Phi(n,k,h)=\Theta(n-1,k,h)$ and define 
\begin{equation}\label{DefEk}
E_h^{k, n} =
\left\{
	\begin{array}{lll}
        0  & \mbox{if } h>k \mbox{ or } h< 0,\\
		1  & \mbox{if } h=k, \\
		-\sum_{i=h+1}^k(-1)^{i-h} \Phi(n,i,h) E_{i}^{k,n} & \mbox{otherwise.}
	\end{array}
\right.
\end{equation}
\begin{proposition}\label{magic5S}
 Let $R_i$ be the reduced recurrence for $C_{2|i}$. Then
 \begin{equation}\label{n=2k+1S5}
 \sum E_i^{k,2k+1}R_i = \Lambda_{2|k}^{2|k, 2k+1} C_{2|k,2k+1}+ \Gamma_{k}^{2|k, 2k+1} C_{k,2k+1}  - \sum_{j=0}^{k-1}s(k-j,q)C_{j,2k+1},
\end{equation}
 \begin{equation}\label{ngenericoS5}
\sum E_i^{k,n}R_i = \Lambda_{2|k}^{2|k, n} C_{2|k,n}+ \Gamma_{k+1}^{2|k, n} C_{k+1,n} + \Gamma_{k}^{k, n} C_{2|k,n}  - \sum_{j=0}^{k-1}d_{k-j,n-2j}C_{j,n}.
 \end{equation}
 \proof
 Coherently with notation of
 Section 4, we will denote by $\Gamma_{2|h}^{2|k,n}$ (resp $\Gamma_{h}^{2|k,n}$) the coefficient of $C_{2|h}$ (resp. $C_h$) in $ \sum E_i^{k,2k+1}R_i$.
 Observing that $E_h^{k,n}=A_h^{k,n-1}$ we obtain an analogous of Lemma \ref{propAhkn}:
 \begin{lemma}\label{propEhkn}
 \begin{enumerate}
  \item $E_h^{k,n}=E_{h-1}^{k-1,n-2}$ $\mathrm{for}$ $h>1$,
  \item $E_{h}^{k,n}=E_{h}^{k,n-1}+E_{h}^{k-1,n-1}$ $\mathrm{if}$ $n \neq 2k+1, 2k+2$,
  \item $E_{h}^{k,2k+1}=E_{h}^{k-1,2k}$,
  \item $E_{h}^{k,2k+2}=2E_{h}^{k,2k+1}+E_{h}^{k-1,2k+1}$.
 \end{enumerate}
\end{lemma}
 An immediate consequence of definition \ref{DefEk} is that $\Gamma_{2|h}^{2|k,n}=0$ if $h<k$. 
If $k>h>0$, we have:
 \begin{align*}
  \Gamma_{h}^{2|k,2k+1}&=\sum_{j=h-1}^{k}E_j^{k,2k+1}\Lambda_h^{2|j,n}\\
  &=\sum_{j=h-1}^{k}E_j^{k,2k+1}\Lambda_{h-1}^{2|j-1,2k-1}-\sum_{j=h-1}^{k-1}E_j^{k,2k+1}\Lambda_{2|h-1}^{2|j,2k}-s(2k-2,q)\left[\sum_{j=h-1}^{k-1}(-1)^{j-h+1}E_j^{k,2k+1}\Phi(2k,j,h-1)\right]\\
  &=\sum_{j=h-1}^{k}E_j^{k,2k+1}\Lambda_{h-1}^{2|j-1,2k-1}-\sum_{j=h-1}^{k-1}E_j^{k-1,2k}\Lambda_{2|h-1}^{2|j,2k}-s(2k-2,q)\left[\sum_{j=h-1}^{k-1}(-1)^{j-h+1}E_j^{k-1,2k}\Phi(2k,j,h-1)\right]\\
  &=\Gamma_{h-1}^{2|k-1,2k-1},
 \end{align*}
\begin{align*}
\Gamma_{h}^{2|k,2k+2}&=\sum_{j=h-1}^{k}E_j^{k,2k+2}\Lambda_{h-1}^{2|j-1,2k}-2\Lambda_{2|h-1}^{2|k,2k+1}-\sum_{j=h-1}^{k-1}E_j^{k,2k+2}\Lambda_{2|h-1}^{2|j,2k+1}+\\&-s(2k-1,q) \left[2(-1)^{k-h+1}\Phi(2k+1,k,h-1) +\sum_{j=h-1}^{k-1}(-1)^{j-h+1}E_j^{k,2k+2}\Phi(2k+1,j,h-1)\right]\\
&=\Gamma_{h-1}^{2|k-1,2k}-\sum_{j=h-1}^{k}\left[2E_j^{k,2k+1}+E_j^{k-1,2k+1}\right]\Lambda_{2|h-1}^{2|j,2k+1}+\\&-s(2k-1,q) \left[\sum_{j=h-1}^k(-1)^{j-h+1}\left(2E_j^{k,2k+1}+E_j^{k-1,2k+1}\right)\Phi(2k+1,j,h-1)\right]\\
&=\Gamma_{h-1}^{2|k-1,2k}-2\sum_{j=h-1}^{k}E_j^{k,2k+1}\Lambda_{2|h-1}^{2|j,2k+1}  -\sum_{j=h-1}^{k-1}E_j^{k-1,n-1}\Lambda_{2|h-1}^{2|j,n-1}\\
&=\Gamma_{h-1}^{2|k-1,2k},
 \end{align*}
 \begin{align*}
\Gamma_{h}^{2|k,n}&=\sum_{j=h-1}^{k}E_j^{k,n}\Lambda_{h-1}^{2|j-1,n-2}-\sum_{j=h-1}^{k}E_j^{k,n}\Lambda_{2|h-1}^{2|j,n-1}-s(n-3,q) \left[\sum_{j=h-1}^k(-1)^{j-h+1}E_j^{k,n}\Phi(n-1,j,h-1)\right]\\
&=\Gamma_{h-1}^{2|k-1,n-2}-\sum_{j=h-1}^{k}\left[E_j^{k,n-1}+E_j^{k-1,n-1}\right]\Lambda_{2|h-1}^{2|j,n-1}\\
&=\Gamma_{h-1}^{2|k-1,n-2}-\sum_{j=h-1}^{k}E_j^{k,n-1}\Lambda_{2|h-1}^{2|j,n-1}  -\sum_{j=h-1}^{k-1}E_j^{k-1,n-1}\Lambda_{2|h-1}^{2|j,n-1}\\
&=\Gamma_{h-1}^{2|k-1,n-2}.
 \end{align*}
We reduced again to the computation of explicit formulae for the zero coefficient in the recurrences.
\begin{align*}
  \Gamma_{0}^{2|k,2k+1}&=\sum_{j=0}^{k}E_j^{k,2k+1}\Lambda_0^{2|j,n}\\
  &=\sum_{j=0}^{k-1}E_j^{k,2k+1}\Lambda_{0}^{2|j,2k}-\sum_{j=1}^{k}E_j^{k,2k+1}\Lambda_{0}^{2|j-1,2k-1}+\\&+\left(s(2k-2,q)+s(-2k,q)\right)\left[\sum_{j=1}^{k}(-1)^{j-1}E_j^{k,2k+1}\Phi(2k+1,j,1)\right]\\
  &=\sum_{j=0}^{k-1}E_j^{k-1,2k}\Lambda_{0}^{2|j,2k}-\Gamma_{0}^{2|k-1,2k-1}\\
  &=\Gamma_{0}^{2|k-1,2k}-\Gamma_{0}^{2|k-1,2k-1},
 \end{align*}
  \begin{align*}
\Gamma_{0}^{2|k,2k+2}&=2\Lambda_0^{2|k,2k+1}+\sum_{j=0}^{k-1}E_j^{k,2k+2}\Lambda_{0}^{2|j,2k+1}-\sum_{j=1}^{k}E_j^{k,2k+2}\Lambda_{0}^{2|j-1,2k}+\\&+\left(s(2k-1,q)+s(-2k-1,q)\right)\left[\sum_{j=1}^{k}(-1)^{j-1}E_j^{k,2k+2}\Phi(2k+2,j,1)\right]\\
  &=2\Lambda_0^{2|k,2k+1}+\sum_{j=0}^{k-1}\left[2E_j^{k,2k+1}+E_j^{k-1,2k+1}\right]\Lambda_{0}^{2|j,2k+1}-\Gamma_{0}^{2|k-1,2k}\\
  &=2\Gamma_{0}^{2|k,2k+1}+\Gamma_{0}^{2|k-1,2k+1}-\Gamma_{0}^{2|k-1,2k},
 \end{align*} 
 \begin{align*}
\Gamma_{0}^{2|k,n}&=\sum_{j=0}^{k}E_j^{k,n}\Lambda_{0}^{2|j,n-1}-\sum_{j=1}^{k}E_j^{k,n}\Lambda_{0}^{2|j-1,n-2}+\\&+\left(s(n-3,q)+s(-n+1,q)\right)\left[\sum_{j=1}^{k}(-1)^{j-1}E_j^{k,n}\Phi(n,j,1)\right]\\
  &=\sum_{j=0}^{k-1}\left[E_j^{k,n-1}+E_j^{k-1,n-1}\right]\Lambda_{0}^{2|j,n-1}-\Gamma_{0}^{2|k-1,n-2}\\
  &=\Gamma_{0}^{2|k,n-1}+\Gamma_{0}^{2|k-1,n-1}-\Gamma_{0}^{2|k-1,n-2}.
 \end{align*} 
A direct computation shows that $\Gamma_0^{2|2,5}=-s(2,q)$. Using the above relations between the $\Gamma_0^{2|k,n}$, the proposition now follows proving by induction the identities
\[\Gamma_0^{2|k,2k+1}=-s(k,q), \qquad \Gamma_0^{2|k,n}=-s(k,q)-s(n-k-1,q).\]
\end{proposition}
\section{The Exceptional Cases}
The most efficient way to verify the Reeder's Conjecture in the Exceptional cases is to compute explicitly the coefficients appearing in Stembridge's recurrences (evaluated in $-q$ and $q^2$) using a computer algebra software (in our case, SageMath 9.3 \cite{Sage}) and prove that the  polynomials computed by Gyoja, Nishiyama and Shimura for \cite{GnS} satisfies such recurrences. 
We underline that a crucial ingredient for these computations is the strong efficiency of the Stembridge's algorithms.

As pointed out in \cite{Stembridge}, for cases $F_4$ and $E_8$ it is not possible to use minuscule recurrences (these algebrae do not have a minuscule coweight) and then we were forced to use Stembridge's \emph{quasi-minuscule recurrence}(\cite{Stembridge}, Formula 5.13), already used to prove Reeder's Conjecture in type $C$(see \cite{SDT}).

We list in the tables below the small weights different from $0$, $\theta_s$ and $\theta$ in the exceptional cases. We used the numbering of Dynkin diagrams as in \cite{Bou}. The description of their zero weight spaces is exposed in \cite{AHJR},\cite{R4} using the notation of \cite{Carter}. Moreover, for the $F_4$ case we used the computations contained in \cite{BL} and the results of Section 4.10 in \cite{Lu}.   
\begin{table}[h!]\label{smallweightsE6}
  \begin{center}
   \caption{ Type $E_6$.  Minuscule coweights: $\omega_1, \omega_6$ }
   \begin{tabular}{c|c}
    \textbf{Small Representation} & \textbf{Zero Weight Space } \\
    Highest weight & $\phi_{\alpha , \beta}$ description \\
    \hline
         $\omega_{1}+\omega_{6}$  &  $\phi_{20,2}$ \\
    $\omega_{4}$ &  $\phi_{30,3}+\phi_{15,5}$ \\
    $\omega_1 + \omega_{3}$, $\omega_5 + \omega_{6}$ & $\phi_{64,4}$\\
    $3\omega_1$, $3\omega_{6}$ & $\phi_{24,6}$
\end{tabular}
\end{center}
\end{table}
\begin{table}[h!]\label{smallweightsE7}
  \begin{center}
   \caption{ Type $E_7$.  Minuscule coweights: $\omega_7$ }
   \begin{tabular}{c|c}
    \textbf{Small Representation} & \textbf{Zero Weight Space } \\
    Highest weight & $\phi_{\alpha , \beta}$ description \\
    \hline
    $\omega_{6}$ &  $\phi_{27,2}$ \\
    $\omega_{3}$ &  $\phi_{56,3}+\phi_{21,6}$ \\
      $2\omega_{7}$  &  $\phi_{21,3}$ \\
    $\omega_2 + \omega_{7}$ & $\phi_{120,4}+\phi_{105,5}$\\
\end{tabular}
\end{center}
\end{table}
\begin{table}[h!]\label{smallweightsE8}
  \begin{center}
   \caption{ Type $E_8$.  Quasi-minuscule coweight: $\omega_8$ }
   \begin{tabular}{c|c}
    \textbf{Small Representation} & \textbf{Zero Weight Space } \\
    Highest weight & $\phi_{\alpha , \beta}$ description \\
    \hline
    $\omega_{1}$ &  $\phi_{35,2}$ \\
      $\omega_{7}$  &  $\phi_{112,3}+\phi_{28,8}$ \\
    $\omega_2$ & $\phi_{210,4}+\phi_{160,7}$\\
\end{tabular}
\end{center}
\end{table}
\begin{table}[h!]\label{smallweightsF4}
  \begin{center}
   \caption{ Type $F_4$.  Quasi-minuscule coweight: $\omega_4$ }
   \begin{tabular}{c|c}
    \textbf{Small Representation} & \textbf{ Zero Weight Space  } \\
    Highest weight & $\phi^{'}_{\alpha , \beta}$ description  \\
    \hline
      $\omega_{3}$ &  $\phi^{'}_{8 , 3}+\phi^{'}_{1 , 12}$ \\
\end{tabular}
\end{center}
\end{table}

We point out that, for type $E_6$, some nice symmetries can be used to reduce the number of cases that need to be examined. We remark that in Table 2 there are two pairs of weights that have the same zero weight space representation. 
 Set $X_1=\{3 \omega_1, \omega_1+\omega_3\}$ and $X_6=\{3 \omega_6, \omega_5+\omega_6\}$. We recall that the automorphism of the $E_6$ Dynkin diagram that send labelled vertices 1 to 6 and 3 to 5 induces an involution $J$ on the weight lattice that pairwise exchange the weights in $X_1$ with weights in $X_2$.
Let us denote by $R_\lambda^{\omega_1}$ (resp. $R_\lambda^{\omega_6}$) the minuscule recurrence for the weight $\lambda$, obtained choosing $\omega_1$ (resp. $\omega_6$) as a minuscule coweight.
Fix $\lambda \in X_1$. We observed empirically that the involution $J$ sends $R_\lambda^{\omega_1}$ to 
 $R_{J(\lambda)}^{\omega_6}$, preserving the coefficients appearing in the recurrences. This proves that $C_\lambda(q)$ and $C_{J(\lambda)}(q)$ satisfy the same recursive relations and then we reduced to prove Reeder's Conjecture only for weights in $X_1$.
Our computations are available at the link in Bibliography \cite{SDTComp}.

\end{document}